\documentclass[reqno]{amsart}
\usepackage{amssymb}
\def\R{{\mathbb R}}
\def\C{{\mathbb C}}
\def\Z{{\mathbb Z}}
\def\V{{\mathbb V}}
\def\a{\alpha}
\def\b{\beta}
\def\sign{\mathop{\rm sign}}
\def\Or{\mathop{\rm Or}}
\def\orr{\mathop{\rm or}}
\def\Aut{\mathop{\rm Aut}}
\newtheorem{thm}{Theorem}[section]
\newtheorem{lem}[thm]{Lemma}
\newtheorem{prop}[thm]{Proposition}
\newtheorem{defin}[thm]{Definition}
\begin{document}
\title{On the orientation of graphs}
\author{Konstantin Salikhov}
\thanks{Partially supported by the Russian Fond for Basic Research
grant 99-01-00009}
\address{Dept. of Differential Geometry, Faculty of Mechanics and
Mathematics, Moscow State University, Moscow, 119899, Russia}
\email{salikhov@mccme.ru}
\maketitle
\section{Introduction}
Soon after the first papers on calculation of (co)homology groups of
Lie algebras appeared in 70s, the idea to use graphs to denote
generators of these groups arose. By the Invariants Theory, the calculation
of these groups can be reduced to the calculation of (co)homology groups of
some graph-complices for a wide class of Lie algebras \cite{F, K}.
The generators of a graph-complex are pairs $(\Gamma,\orr_\Gamma)$,
where $\Gamma$ is a connected graph, and $\orr_\Gamma$ is an element of
the $\Z_2$-module $\Or(\Gamma)$ of {\it orientations of $\Gamma$}.
The definition of $\Z_2$-module $\Or(\Gamma)$ contains all the 
information on the (co)homologies of the graph-complex.
The first example of a graph-complex appeared in \cite{K}. It
was constructed to calculate homology groups of the algebra $h_\infty$
($h_\infty$ is the Lie algebra of Hamiltonian vector
fields in infinite-dimensional even real vector space, vanishing at the
origin). In \cite{KS} another graph-complex (with a new definition of
$\Z_2$-module $\Or(\Gamma)$) was constructed to calculate the
(co)homology of the same Lie algebra, and it was conjectured that this
definition coincided with that from \cite{K}.

In this short notice we give a universal definition of $\Z_2$-module
$\Or(\Gamma)$ of orientations of a graph $\Gamma$ and construct a method,
by means of which one can easily verify whenever two such special definitions
coincide. In particular, we prove the conjecture from \cite{KS}.
\section{Definitions and formulation of the result}
To describe graphs we will use the so-called {\it half-edge language}.
A graph $\Gamma$ is a set of {\it half-edges} $H(\Gamma)$, consisting of $2n$ elements,
and two partitions of this set into disjoint unions of subsets. The first
partition $E(\Gamma)$ is called the set of {\it edges} and consists of $n$
two-element subsets of $H(\Gamma)$, and the second partition $V(\Gamma)$ is
called the set of {\it vertices} and consists of arbitrary non-empty subsets.
An edge is called a {\it loop} if both its half-edges belong to one vertex.
An {\it automorphism} of the graph $\Gamma$ is an arbitrary permutation on the
set $H(\Gamma)$, which respects the structures of $E(\Gamma)$ and $V(\Gamma)$.

Let $\Gamma$ be a graph and $\Aut(\Gamma)$ be its group of automorphisms.
Let us call by an {\it orientation homomorphism} an arbitrary homomorphism
$\Theta:\Aut(\Gamma)\to\Z_2$ (here and in the sequel $\Z_2=\{+1,-1\}$).
Consider the set of all pairs $(\tau,\varepsilon)$, where $\tau$ is
an enumeration of the elements of $V(\Gamma)$ by the numbers
$1,\dots,|V(\Gamma)|$, and $\varepsilon\in\Z_2$.
The group $\Aut(\Gamma)$ acts on this set by the formula
\begin{equation*}
\varphi(\tau,\varepsilon)=
(\varphi(\tau),\Theta(\varphi)\cdot\varepsilon),
\end{equation*}
where $\varphi\in\Aut(\Gamma)$, $\varphi(\tau)$ is the induced enumeration
of $V(\Gamma)$ and $\Theta(\varphi)$ is the value of the orientation
homomorphism $\Theta:\Aut(\Gamma)\to\Z_2$ on the element
$\varphi\in\Aut(\Gamma)$. If this action is degenerated, the graph
$\Gamma$ is called {\it $\Theta$-non-orientable}, and $\Or(\Gamma)=0$.
If this action is non-degenerated, the graph $\Gamma$ is called
{\it $\Theta$-orientable}, and $\Or(\Gamma)$ is the set of left
cosets of this action. The action of $\Z_2$ on $\Or(\Gamma)$
is defined by the formula
$\delta(\tau,\varepsilon)=(\tau,\delta\cdot\varepsilon)$.
Note that if the graph $\Gamma$ is 1-skeleton of $n$-simplex
and $\Theta:\Aut(\Gamma)\to\Z_2$ is the parity of the induced vertex
permutation, this definition becomes the usual definition of orientation
of $n$-simplex.
\begin{defin}[\cite{K}]\label{defKon}
Let $\Gamma$ be a graph. An automorphism $\varphi\in\Aut(\Gamma)$ induce a
permutation $\pi_\varphi$ of edges of the graph $\Gamma$ and a
non-degenerated linear transformation
$A_\varphi\in\mathop{\rm GL}(H_1(\Gamma,\R))$ in the first homology group
of the graph $\Gamma$ (as a topological space) with real coefficients.
The Kontsevich orientation homomorphism $\Theta_K:\Aut(\Gamma)\to\Z_2$
is defined by the formula
$\Theta_K(\varphi)=\sign(\pi_\varphi)\cdot\sign(\det A_\varphi)$.
\end{defin}
\begin{defin}[\cite{KS}]\label{defKS}
Let $\Gamma$ be a graph. Let us arrange arrows on the edges of $\Gamma$
arbitrarily. An automorphism $\varphi\in\Aut(\Gamma)$ induces
a permutation $\sigma_\varphi$ of vertices of the graph  $\Gamma$
and a function $\varepsilon_\varphi:E(\Gamma)\to\Z_2$
(for each edge $e\in E(\Gamma)$ the value $\varepsilon_\varphi(e)=+1$,
if the orientation of the arrow on $e$ coincides with the induced
orientation $\varphi_*(\varphi^{-1}(e))$, and $\varepsilon_\varphi(e)=-1$
otherwise). The Shoikhet orientation homomorphism
$\Theta_S:\Aut(\Gamma)\to\Z_2$ is defined by the formula
\begin{equation*}
\Theta_S(\varphi)=\sign(\sigma_\varphi)\cdot
\prod_{e\in E(\Gamma)}\varepsilon_\varphi(e)
\end{equation*}
It is easy to see that this product is independent on the arrangement
of arrows on the edges $E(\Gamma)$.
\end{defin}
Note that any graph, containing a loop, is $\Theta_S$- and $\Theta_K$-
non-orientable. To prove this, it suffices to consider the automorphism
of the graph $\Gamma$, which interchanges these two half-edges.
\begin{thm}\label{result}
Let $\Gamma$ be a graph. Then for any $\varphi\in\Aut(\Gamma)$
we have $\Theta_K(\varphi)=\Theta_S(\varphi)$.
\end{thm}
\section{Proof of theorem~\ref{result}}
\begin{proof}[Sketch of the proof of theorem~\ref{result}]
An automorphism $\varphi\in\Aut(\Gamma)$ induces a permutation of edges
of the graph  $\Gamma$. This permutation splits into several cycles.
Let $e$ be an edge of $\Gamma$ and
$\bigl(e,\varphi(e),\dots,\varphi^N(e)\bigr)$ be one of these cycles.
Consider the graph $\Gamma'=\Gamma/e/\varphi(e)/\dots/\varphi^N(e)$,
which is obtained from $\Gamma$ by a contraction of edges
$e,\varphi(e),\dots,\varphi^N(e)$ (we assume that the contraction of
a loop is just a deletion).
The automorphism $\varphi\in\Aut(\Gamma)$ induces the automorphism
$\varphi'\in\Aut(\Gamma')$. By the induction hypothesis,
$\Theta_K(\varphi')=\Theta_S(\varphi')$. Hence, to prove theorem~\ref{result},
it suffices to show that
\begin{equation}\label{main}
\frac{\Theta_K(\varphi)}{\Theta_K(\varphi')}=
\frac{\Theta_S(\varphi)}{\Theta_S(\varphi')}
\end{equation}
The edges $e,\varphi(e),\dots,\varphi^N(e)$ form a subgraph $G\subset\Gamma$.
We construct three families of relatively simple graphs and show that
we have to prove the identity~(\ref{main}) only for the case when
$G$ is a graph from these families.
\end{proof}
\begin{prop}\label{power2}
An automorphism $\varphi\in\Aut(\Gamma)$ induces permutations of vertices,
edges and half-edges of the graph  $\Gamma$. All these permutations split
into several cycles. Then, without loss of generality we may assume that
the lengths of all these cycles are powers of the number 2.
\end{prop}
\begin{proof}
Since $\Theta:\Aut(\Gamma)\to\Z_2$ is a homomorphism, for any odd
$N$ and any $\varphi\in\Aut(\Gamma)$ we have
$\Theta(\varphi^N)=\Theta(\varphi)$. If $N$ is divisible by all
odd factors of the lengths of all these cycles, then the lengths of all
cycles for the automorphism $\varphi^N$ will be powers of 2.
\end{proof}
Choose any cycle from the cycles into which the permutation of edges splits.
From proposition~\ref{power2} it follows that the length of this cycle
is~$2^n$ for some integer $n$. The edges from this cycle form a subgraph
$G$ of the graph $\Gamma$. Note that the automorphism $\psi=\varphi|_G$
generates a cyclic subgroup in $\Aut(G)$, which acts on
the set $E(G)$ transitively.
\begin{lem}
Let $G$ be a graph with $2^n$ edges and an automorphism $\psi\in\Aut(G)$
generate a cyclic subgroup $\Psi=\{1,\psi,\psi^2,\dots\}$ in $\Aut(G)$,
which acts on the set of edges $E(G)$ of the graph $G$ transitively.
Then the pair $(G,\psi)$ coincides with one of the following. Here
$\a$ and $\b$ are symbols to denote half-edges; $n$ (the "number" of edges),
$c$ (the "number" of connected components) and $m$ (the "mass" of each
connected component) are fixed integers; and $i,j,k,l$ are indices.
\begin{itemize}
\item[($i$)]
\ \par\noindent
\begin{tabular}{l}
$H(G)=\bigcup_{i,j}\a_i^j$,
for $i=1\dots2^{n-c+1}$, $j=1\dots2^c$;\\
$E(G)=\bigsqcup_{i,j}\{\a_i^j,\a_{i+2^{n-c}}^j\}$,
for $i=1\dots2^{n-c+1}$, $j=1\dots2^c$;\\
$V(G)=\bigsqcup_j\{\bigcup_i \a_i^j\}$,
for $i=1\dots2^{n-c+1}$, $j=1\dots2^c$;\\
$\psi(\a_i^j)=\a_i^{j+1}$, for $i=1\dots2^{n-c+1}$, $j=1\dots2^c-1$;\\
$\psi(\a_i^{2^c})=\a_{i+1}^1$, for $i=1\dots2^{n-c+1}-1$;
$\psi(\a_{2^{n-c+1}}^{2^c})=\a_1^1$;\\
\end{tabular}
\item[($ii$)]
\ \par\noindent
\begin{tabular}{l}
$H(G)=\bigcup_{i,j}\a_i^j\cup\bigcup_{i,j}\b_i^j$,
for $i=1\dots2^{n-c}$, $j=1\dots2^c$;\\
$E(G)=\bigsqcup_{i,j}\{\a_i^j,\b_i^j\}$,
for $i=1\dots2^{n-c}$, $j=1\dots2^c$;\\
$V(G)=\bigsqcup_j\{\bigcup_i \a_i^j\}\sqcup
\bigsqcup_{j,k}\{\bigcup_l \b_{k+l\cdot 2^m}^j\}$,
for $i=1\dots2^{n-c+1}$, $j=1\dots2^c$,\\
\hskip 15mm $c+m\le n$, $k=1\dots 2^m$, $l=1\dots2^{n-c-m}$;\\
$\psi(\a_i^j)=\a_i^{j+1}$ and $\psi(\b_i^j)=\b_i^{j+1}$,
for $i=1\dots2^{n-c}$, $j=1\dots2^c-1$;\\
$\psi(\a_i^{2^c})=\a_{i+1}^1$ and $\psi(\b_i^{2^c})=\b_{i+1}^1$,
for $i=1\dots2^{n-c}-1$;\\
$\psi(\a_{2^{n-c}}^{2^c})=\a_i^1$ and $\psi(\b_{2^{n-c}}^{2^c})=\b_1^1$;\\
\end{tabular}
\item[($iii$)]
\ \par\noindent
\begin{tabular}{l}
$H(G)=\bigcup_{i,j}\a_i^j\cup\bigcup_{i,j}\b_i^j$,
for $i=1\dots2^{n-c}$, $j=1\dots2^c$;\\
$E(G)=\bigsqcup_{i,j}\{\a_i^j,\b_i^j\}$,
for $i=1\dots2^{n-c}$, $j=1\dots2^c$;\\
$V(G)=\bigsqcup_{j,k}\{\bigcup_l \a_{k+l\cdot 2^m}^j\cup
\bigcup_l \b_{k-1+l\cdot 2^m}^j\}$,
for $i=1\dots2^{n-c+1}$, $j=1\dots2^c$,\\
\hskip 15mm $c+m\le n$, $k=1\dots 2^m$, $l=1\dots2^{n-c-m}$;\\
$\psi(\a_i^j)=\a_i^{j+1}$ and $\psi(\b_i^j)=\b_i^{j+1}$,
for $i=1\dots2^{n-c}$, $j=1\dots2^c-1$;\\
$\psi(\a_i^{2^c})=\a_{i+1}^1$ and $\psi(\b_i^{2^c})=\b_{i+1}^1$,
for $i=1\dots2^{n-c}-1$;\\
$\psi(\a_{2^{n-c}}^{2^c})=\a_i^1$ and $\psi(\b_{2^{n-c}}^{2^c})=\b_1^1$;\\
\end{tabular}
\end{itemize}
\end{lem}
\begin{proof}
Suppose the graph $G$ has a loop. Since the action of $\Psi$ is transitive
on $E(G)$, all the edges of $G$ are loops. Note that there exist not more
than two equivalence classes of half-edges $H(G)$ with respect to the action
of $\Psi$. If this action is transitive, we have the case~$(i)$. If this action
is non-transitive, we have the case~$(iii)$ for $m=0$. In the sequel, we will
assume that the graph $G$ has no loops.

To each of $2^n$ edges of the graph $G$ we associate a point of the kind
$e^{2\pi ip/2^n}$ on the circle $|z|=1, z\in\C$ in such a way that
the action $\psi$ on $E(G)$ corresponds to the multiplication of
the points on the circle by $e^{2\pi i/2^n}$. For each $v\in V(G)$,
consider the subset $E_v(G)\subset E(G)$ of those edges of the graph $G$
that have at least one half-edge in common with the vertex $v$, and consider
the polygon $P_v$ with vertices in the points of the kind $e^{2\pi ip/2^n}$,
corresponding to the elements of the subset $E_v(G)$. Even if for some
$v_1\neq v_2\in V(G)$ the polygons $P_{v_1}$ and $P_{v_2}$ coincide
as the subsets of $\C$, nevertheless we will consider them distinct.
Denote the set of such polygons by $\V$.
Since any edge contains two half-edges and the graph $G$ has no loops,
any point of the kind $e^{2\pi ip/2^n}$ is a vertex of exactly two polygons
from $\V$. Since $\psi\in\Aut(G)$, for any $P\in\V$ we have
$e^{2\pi i/2^n}(P)\in\V$.

Since the action of $\Psi$ on $E(G)$ is transitive, there exist not more
than two equivalence classes of vertices $V(G)$ with respect to the action
of $\Psi$. Suppose there exist two equivalence classes of vertices $V(G)$
with respect to this action. Then the set of polygons $\V$ splits into two
subsets $\V_1\sqcup\V_2$ in such a way that the action of $\Psi$ on each of
$\V_1$ and $\V_2$ is transitive, and any point of the kind $e^{2\pi ip/2^n}$
belongs to exactly one polygon from $\V_1$ and exactly one polygon from $\V_2$.
Since $e^{2\pi i/2^n}(P)\in\V$ for any $P\in\V$, both families $\V_1$ and
$\V_2$ consist of right polygons, and we obtain the case~$(ii)$.

Suppose the action of $\Psi$ on $V(G)$ is non-transitive. Then all
polygons from $\V$ are congruent. Let us call by length of a side of
a polygon the central angle, based on this side, multiplied by $2^n/{2\pi}$.
Note that there exist not more than two different lengths of sides
among all lengths of the sides of all polygons from $\V$.
If all these lengths are equal,
we obtain the case~$(iii)$ for $m=0$. If there exist two different
lengths of sides of polygons from $\V$, we obtain the case~$(iii)$ for $m>0$.
\end{proof}
\begin{proof}[Proof of theorem~\ref{result}]
The identity~(\ref{main}) is proved by straightforward calculations.
To illustrate the idea of these calculations, we consider here only
the case~$(i)$ for $c=0$. All other cases are proved analogously.
\par\noindent
\begin{enumerate}
\item The permutation of edges $\pi_\varphi$ (cf. definition~\ref{defKon})
splits into sum of permutations $\pi_{\varphi'}$ and $\pi_\psi$. Therefore
$\sign(\pi_\varphi)/\sign(\pi_{\varphi'})=\sign(\pi_\psi)=-(-1)^{2^n}$.
\item The permutation of vertices $\sigma_\varphi$ coincides with
$\sigma_{\varphi'}$ (cf. definition~\ref{defKS}). Therefore
$\sign(\sigma_\varphi)=\sign(\sigma_{\varphi'})$.
\item In the group $H_1(G,\R)$, consider the basis
$\{a_i=(\a_i,\b_i)\mid i=1\dots2^n\}$. The inclusion
$G\hookrightarrow\Gamma$ induces the monomorphism
$h:H_1(G,\R)\to H_1(\Gamma,\R)$. Let us complete the set
$\{h(a_i)\mid i=1\dots2^n\}$ to a basis in $H_1(\Gamma,\R)$. In this
basis $A_\varphi=A_{\varphi'}\oplus A_\psi$. Therefore
\begin{equation*}
\sign(\det A_\varphi)/\sign(\det A_{\varphi'})=\sign(\det A_\psi)=
\left\|
\begin{array}{ccccc}
0 & 0 & \cdots & 0 & 1\\
1 & 0 & \cdots & 0 & 0\\
0 & 1 & \cdots & 0 & 0\\
\vdots & \vdots & \ddots & \vdots & \vdots \\
0 & 0 & \cdots & 1 & 0
\end{array}\right\|=-(-1)^{2^n}
\end{equation*}
\item On the edges $\{\a_i,\b_i\}$ of the graph $G$, arrange the arrows
$\a_i\to \b_i$ for $i=1\dots2^n$. Then for any edge $e\in E(G)$ we have
$\varepsilon_\psi(e)=+1$. Therefore
\begin{equation*}
\prod_{e\in E(\Gamma)}\varepsilon_\varphi(e)\bigg/
\prod_{e\in E(\Gamma')}\varepsilon_{\varphi'}(e)=
\prod_{e\in E(G)}\varepsilon_\psi(e)=1
\end{equation*}
\end{enumerate}
From items $1\dots4$ it follows that
${\Theta_K(\varphi)}/{\Theta_K(\varphi')}=
{\Theta_S(\varphi)}/{\Theta_S(\varphi')}$.
\end{proof}
\section*{Acknowledgments}
I would like to thank S.V.~Duzhin for suggesting the problem and
B.~Shoikhet for useful discussions.

\end{document}